\documentclass[a4paper,11pt]{article}

\usepackage[margin=1in]{geometry}  

\usepackage{moreverb}

\usepackage[dvips,colorlinks,bookmarksopen,bookmarksnumbered,citecolor=red,urlcolor=red]{hyperref}

\usepackage{amsfonts,amsmath,amssymb,amsthm}
\usepackage{graphicx}

\usepackage{algorithm}
\usepackage{algorithmic}

\newcommand{\NN}{\mathbb{N}}
\newcommand{\RR}{\mathbb{R}}

\begin{document}

\title{Hermite-Birkhoff Interpolation on Arbitrarily Distributed Data on the Sphere and Other Manifolds}

\author{Giampietro Allasia, Roberto Cavoretto, Alessandra De Rossi\footnote{Department of Mathematics \lq\lq G. Peano\rq\rq, University of Torino, via Carlo Alberto 10, I--10123 Torino, Italy. E-mails: giampietro.allasia@unito.it, roberto.cavoretto@unito.it, alessandra.derossi@unito.it}}

\date{}

\maketitle

\begin{abstract}
We consider the problem of interpolating a function given on scattered points using Hermite-Birkhoff formulas on the sphere and other manifolds. We express each proposed interpolant as a linear combination of basis functions, the combination coefficients being incomplete Taylor expansions of the interpolated function at the interpolation points. The basis functions have the following features: (i) depend on the geodesic distance; (ii) are orthonormal with respect to the point-evaluation functionals; and (iii) have all derivatives equal zero up to a certain order at the interpolation points. Moreover, the construction of such interpolants, which belong to the class of partition of unity methods, takes advantage of not requiring any solution of linear systems. 
\end{abstract}





\section{INTRODUCTION}

In previous papers we dealt with Hermite-Birkhoff interpolation of a function given on arbitrarily distributed points on Euclidean spaces \cite{bracco11b} (see e.g. also \cite{caira12, costabile12, costabile13, dellaccio16}), and with Lagrange and Hermite-Birkhoff interpolation on Banach spaces \cite{bracco11c}. The Hermite and generalized Hermite  interpolations on scattered data by means of basis functions depending on the distance have been proposed by Wu's pioneering paper in 1992 \cite{wu92}. Since then, the interest in this topic seems to have increased significantly (see  e.g. the pertinent chapters  in \cite{fass07, wendl05}). A number of authors have also considered the  Hermite  interpolation setting on scattered data  on the sphere (see e.g. \cite{ron96, fass99, narc02, macedo11}) or even general Riemannian manifolds \cite{dyn99, narc95}. The sphere is a particularly interesting example of a connected compact smooth manifold, even because considerations developed about interpolation on the sphere can be extended to other manifolds. We think convenient to discuss, as long as possible, a general framework, though in practice  the most interesting manifolds are smooth two-dimensional manifolds, i.e. surfaces in $\RR^3$. 

Here we consider the cardinal (radial) basis function method. In this way, in the Hermite-Birkhoff interpolation the interpolant is directly expressed  as a linear combination of  basis functions, which depend on the geodesic distance, are orthonormal with respect to the point-evaluation functionals, and have all  derivatives equal zero up to a certain order at the interpolation points. The  coefficients of the linear combination are incomplete Taylor expansions of the interpolated function at the interpolation points.  
Our interpolation method is strictly linked up with papers which discuss Lagrange interpolation  by  partition of unity methods, namely Shepard-like methods, on the sphere (see, in particular, \cite{cavoretto10b, cavoretto11}).
In this paper explicit expressions of Hermite-Birkhoff interpolants on manifolds are given. Since these definitions are based on a suitable class of cardinal basis functions, we provide a general way to construct such basis functions, which depend on geodesic distances on Riemannian manifolds. Moreover, upper bounds for errors in terms of the fill distance are shown. 


\section{HERMITE-BIRKHOFF INTERPOLATION ON MANIFOLDS} 


\noindent
{\bf Definition 1.} {\sl
Let us consider  a  $m$-dimensional Riemannian manifold ${\cal M}\subset\RR^{d+1}, ~d\ge 2$,  simply connected, compact, 
an open set $U=\{u\equiv (u_1, \ldots, u_{d+1})\in\RR^{d+1}\}\subset {\cal M}$, a function $\varphi:U\to\RR^d$ which maps $U$ homeomorphically to the open set $V=\{v\equiv (v_1, \ldots, v_m)\in\RR^m\}:=\varphi (U)$ so that $v=\varphi (u):=v(u)$ and $u=\varphi^{-1}(v):=u(v)$.
Let ${\cal X}=\lbrace  z_1,..., z_n\rbrace\subset U$ be a set of  distinct points, possibly scattered,  with associated finite sets $\Delta_1,...,\Delta_n\subset \NN_0^m$, always containing the value 0. The {\rm Hermite-Birkhoff interpolation problem} from  $U$ to $\RR$  consists in finding a function  $H:U\rightarrow \RR$ which satisfies the interpolation conditions
\begin{equation}
D^\beta H( z_i):=\frac{\partial^{|\beta|}H( z_i)}{\partial v_1^{\beta_1}\cdots \partial v_m^{\beta_m}}=f_{i \beta}, \qquad \beta\in\Delta_i, 
\quad i=1,...,n,
\end{equation}
where $\beta=(\beta_1, \ldots, \beta_m)$, $\vert\beta\vert=\beta_1+\cdots +\beta_m$, and the $f_{i \beta}\in \RR$ are  given values to be interpolated.
It is  assumed that $H\in C^k(U)$, where $k=\max\lbrace\vert \beta\vert: \beta\in\Delta_i,\,\hbox{for some }i,\, 1\le i\le n\rbrace$.
}
\vskip 0.1cm

  In the following, we will also think that the $f_{i\beta}$ are values assumed by an underlying  function $f:U\rightarrow\RR$, $f\in{\cal C}^k(U)$, so that the conditions (1)  take the form
\begin{equation}
D^\beta H( z_i)=D^\beta f( z_i), \qquad \beta\in\Delta_i, \quad i=1,...,n.
\end{equation}
In general, the values of $f$ and of some of its derivatives are known only at the points of ${\cal X}$.

A constructive solution to the interpolation problem (1) can be given by introducing a suitable class of {\sl cardinal basis functions}, which can be   defined as follows.
\vskip 0.1cm

\noindent
{\bf Definition 2.} {\sl
\noindent Given a set of distinct points ${\cal X}=\{ z_i,~1\le i\le n \}$, arbitrarily distributed in the open set  $U\subset {\cal M}$, the functions $g_j: U\to\RR,~1\le j\le n$, are cardinal basis functions with respect to ${\cal X}$ if they satisfy for all $u\in U$ the conditions
$$
g_j\in C^k(U),\qquad g_j(u)\ge 0  ,\qquad \sum_{j=1}^n g_j(u)=1 ,\qquad g_j( z_i)=\delta_{ji} ,
$$
where  
$\delta_{ji}$ is the Kronecker delta, and  the additional property
\begin{equation}
D^\beta g_j(z_i)=0, \qquad \beta\in\Delta_i,~\vert \beta\vert\ne 0, \quad i=1,...,n. 
\end{equation}
}
\vskip 0.1cm

It is clear that an interpolant based on these weights must be considered as a {\sl partition of unity method}.
\vskip 0.1cm

\noindent
{\bf Property 3.} {\sl
The interpolation conditions (2) are satisfied  by the interpolant
\begin{equation}
H(u)=\sum_{i=1}^n ~T(u;f, z_i,\Delta_i)~g_i(u),
\end{equation}
where
\begin{eqnarray}
T(u;f,  z_i,\Delta_i):=\sum_{\beta\in\Delta_i} \frac{D^\beta f(z_i)}{\beta_1!\cdots\beta_m!}  \big(v-v(z_i)\big)^\beta
\nonumber
\end{eqnarray}
is formally an incomplete Taylor expansion of $f$ at $ z_i$, in the sense that it only includes the partial derivatives whose orders belong to $\Delta_i$. The interpolant (4) can also be be expressed in the form
\begin{equation}
H(u)=\sum_{i=1}^n \sum_{\beta\in\Delta_i} ~D^\beta f(z_i)~g_{i\beta}(u),
\end{equation}
where
\begin{equation}
g_{i\beta}(u):=
\frac{\big(v-v(z_i)\big)^\beta}{\beta_1!\cdots\beta_m!}g_i(u)=
\frac{\big(v_1-v_1(z_i)\big)^{\beta_1}\cdots \big(v_m-v_m(z_i)\big)^{\beta_m}}{\beta_1!\cdots\beta_m!}g_i(u),  
\end{equation}
with $\beta\in\Delta_i,~ 1\le i\le n$.
}
\vskip 0.1cm

The formula (5) highlights that the interpolant is essentially constructed by using the $g_{i\beta}$  as basis functions.
Moreover, the interpolant (4)
enjoys the usual properties of cardinal basis interpolants.
\vskip 0.1cm

\noindent
{\bf Property 4.} {\sl
\noindent There hold  the following inequalities:
\begin{enumerate}
\item[a)] $\Vert H(u)\Vert\le \max_i \Vert T(u;f, z_i,\Delta_i)\Vert$, \vskip.1cm
\item[b)] $\Vert f(u)-H(u)\Vert\le \sum_{i=1}^n g_i(u)\Vert f(u)-T(u;f,  z_i,\Delta_i)\Vert\le \max_i \Vert f( u)-T( u;f,  z_i,\Delta_i)\Vert$,
\end{enumerate}
where the $i$-th term in the sum may be interpreted as the local error at the point $ z_i$.
}


\section{CARDINAL BASIS FUNCTIONS ON MANIFOLDS}

A classical way to construct cardinal basis functions defined on $\RR^{d+1}$ is Cheney's method (see \cite{cheney} and \cite{cheney00}, pp. 67-68), which can be used for manifolds as well, if we adopt a suitable  distance.  
\vskip 0.1cm

\noindent
{\bf Theorem 5.} {\sl
Let us consider $U\subset {\cal M}$ as in Definition 1 and let $\alpha:U\times U \rightarrow \RR^+$ be a continuous and bounded function, such that $\alpha(u,z_i)>0$ for all $u\in U$, $u\ne z_i$, and $\alpha (z_i, z_i)=0$ for all $z_i\in {\cal X}$. Moreover, let each $\alpha(u,z_i)$ be $k$-times continuously differentiable on $U$ such that
\begin{eqnarray*}
[D^{\beta} \alpha(u,z_i)]_{u=z_i}=0, \qquad i=1,2,\ldots,n, \qquad 0<\vert\beta\vert\le k. 
\end{eqnarray*}
The corresponding {\rm cardinal basis functions}
\begin{equation}
g_i(u)={{\displaystyle{\prod_{j=1, j\neq i}^n \alpha (u,z_j)}}
\over{{\displaystyle\sum_{k=1}^n\prod_{j=1, j\neq k}^n \alpha (u,z_j)}}},\qquad i=1,2,\ldots,n,  
\label{ca}
\end{equation}
are continuous and satisfy
\begin{equation}
D^\beta g_i(z_j)=0, \qquad 0<\vert\beta\vert\le k, \quad i,j=1,\ldots,n. 
\label{cb}
\end{equation}
}

\smallskip
\noindent
{\bf Proof:} This result is essentially the $d-$dimensional case of the main theorem in \cite{bracco11b}. 
\hfill $\square$
\vskip 0.1cm


A natural choice is defining $\alpha$ using the distance between points. Since we are considering points on the manifold ${\cal M}$, we   take  the geodesic distance $d_g$ 
and  define  $\alpha$ in the general form
\begin{equation}
\alpha(u, w)=\vartheta (d_g(u, w))  ,
\label{cd} 
\end{equation}
which obviously must satisfy the assumptions of Theorem 5.  In particular, in (9) we  may consider the choice
\begin{equation}
\alpha(u,w)= (d_g(u,w))^\mu, \qquad  \mu\in\RR^+ ,~ \mu\ge k,  ~u,w\in U,
\end{equation}
 which ensures both the vanishing of the derivatives at the nodes and the regularity assumptions, and among the  possible choices is the most direct. Other interesting choices can be found in \cite{cavoretto10a}, where exponential weights that are rapidly decreasing are there considered.

As a result  of the choice (10), we obtain the cardinal basis functions
\begin{equation}
g_i(u)=\frac{(d_g(u, z_i))^{-\mu}}{\sum_{k=1}^n (d(u, z_k))^{-\mu}}, \quad i=i, \ldots, n.
\label{cf}
\end{equation}

For computational reasons, in many cases it may be  preferable to use a localized version of the cardinal basis functions (11), i.e.,
\begin{equation}
\tilde g_i(u)=\frac{\tau_i(u)(d_g(u, z_i))^{-\mu}}{\sum_{k=1}^n \tau_k(u)(d_g(u, z_k))^{-\mu}}, \label{cg}
\end{equation}
where $\tau_i:U\to\RR^+$, $\tau_i\in{ C}^k(U)$, such that
\begin{equation}
\tau_i(u)=
\left\{
\begin{array}{ll}
>0,& ~\hbox{for }  u:d_g(u,z_i)<\delta,\\ 
=0,& ~\hbox{for }  u:d_g(u,z_i)\ge\delta, \\
\end{array}
\right.
\label{cga} 
\end{equation}
and  $\delta>0$ is a suitably chosen value. 

For the  Hermite-Birkhoff interpolant with cardinal basis functions (12)
\begin{equation}
\tilde H(u)=\sum_{i=1}^n ~T(u;f,z_i,\Delta_i)~ \tilde g_i(u),
\end{equation}
we can give more significant error estimates than for the basic case (4).
Let
$q\in\NN $ be defined such that each Taylor-type expansion $T(u; f, z_j, \Delta_j)$ is a complete Taylor expansion up to order $q$, plus other terms of higher degree.
For any $f: U\to\RR$ with $f\in{ C}^q(U)$ and for any $u\in U$, we have, since the cardinal basis functions $\tilde g_i$ are a partition of unity,
\begin{eqnarray}
&\vert f(u)-\tilde H(u)\vert=\Bigg\vert\displaystyle{\sum_{i=1}^n~f(u)\,\tilde g_{i}(u)-\sum_{i=1}^n~T(u; f, z_i, \Delta_i)\,\,\tilde g_{i}(u)}\Bigg\vert=\nonumber\\ 
&\Bigg\vert \displaystyle{\sum_{i=1}^n~\big[f(u)-T(u; f, z_i, \Delta_i)\big]\,\,\tilde g_{i}(u)}\Bigg\vert\le \displaystyle{\sum_{i=1}^n~\big\vert f(u)-T(u; f, z_i, \Delta_i)\big\vert\,\,\tilde g_{i}(u)}, 
\end{eqnarray}
each $\tilde g_i$ being  non-zero only inside the ball of radius $\delta$ centered at $z_i$. Now, since each $T(u;f,z_i,\Delta_i)$ is a Taylor expansion complete up to order  $q$, we can use the estimate 
\begin{equation}
\vert f(u)-T(u; f, z_i, \Delta_i)\vert \le C_i \Vert v-v(z_i)\Vert^{q+1}, 
\end{equation}
where $C_i\in\RR^+$ is a suitable constant and  $\Vert\cdot\Vert$ is the Euclidean norm.
Since $\Vert v-v(z_i)\Vert$ is less than or equal to the geodesic distance $d_g(u, z_i)$, 
it  follows
\begin{equation}
\vert f(u)-T(u; f, z_i, \Delta_i)\vert \le C_i d_g^{q+1}(u,z_i). 
\end{equation}
Inserting (17) in (15) and exploiting again the partition of unity property, since $\vert d_g(u,z_i)\vert<\delta$, we obtain 
\begin{eqnarray*}
\vert f(u)-\tilde H(u)\vert\le \sum_{i=1}^n C_i d_g^{q+1}(u,z_i)\tilde g_i(u)
\le \delta^{q+1}\sum_{i=1}^n C_i \tilde g_i(u)\le C\delta^{q+1}, 
\end{eqnarray*}
with $ C=\max_i C_i$.
Moreover, if we set the localization radius $\delta=K h_{U,{\cal X}}$, where $h_{U,{\cal X}}$ is the so-called {\sl fill distance}, that is,
\begin{equation}
h_{U,{\cal X}}:=\sup_{u\in U }\inf_{z_i\in {\cal X}} d_g(u,z_i).
\end{equation}
and $K\ge 1$, we get the estimate
\begin{equation}
\vert f(u)-\tilde H(u)\vert\le CK h_{U,{\cal X}}^{q+1}. 
\end{equation}


\section{ACKNOWLEDGMENTS}
This research was partially supported by the project ``Metodi numerici nelle scienze applicate'' of the Department of Mathematics \lq\lq Giuseppe Peano\rq\rq\ of the University of Torino, 2014.



\nocite{*}
\bibliographystyle{elsart-num-sort}%
\bibliography{sample_ad}%
\end{document}